\documentclass[10pt]{amsart}

\usepackage{amsmath}
\usepackage{amssymb}
\usepackage{amsthm}
\usepackage {amscd}
\usepackage{epsfig}
\newtheorem{thm}{Theorem}
\newtheorem{defn}{Definition}
\newtheorem{cor}{Corollary}

\theoremstyle{definition}
\newtheorem{example}{Example}

\theoremstyle{remark}
\newtheorem{remark}{Remark}

\newcommand{\C}{{\mathbb C}}

\newcommand{\Z}{{\mathbb Z}}

\newcommand{\CP}{\C P}
\newcommand{\cpq}{\overline{\C P^2}}

\title{On the kinkiness of closed braids}
\author{Christian Bohr}
\address{Mathematisches Institut \\ Theresienstr. 39 \\ 80333 M\"unchen \\ Germany}
\email{bohr@rz.mathematik.uni-muenchen.de}
\begin{document}
\date{\today}
\begin{abstract}
In this note, we prove a lower bound for the positive kinkiness of a closed braid
which we then use to derive an estimate for the positive kinkiness of a link in terms
of its Seifert system. As an application,  we show that certain pretzel knots cannot be unknotted using only 
positive crossing changes. We also describe a subgroup of infinite rank in the smooth knot concordance group 
of which no element has a strongly quasipositive representative.
\end{abstract}
\keywords{Unknotting numbers, kinkiness, strongly quasipositive knots}

\maketitle

In 1993, L. Rudolph proved a lower bound for the slice genus of a knot in terms of a
presentation as the closure of a braid. It is clear that the same estimate holds for the number
of double points of any properly immersed disk in the 4--ball spanning the knot,
for such a disk which has $r$ self--intersection points can be turned into an embedded 
surface of genus $r$ by replacing all the self--intersection points by handles. In this paper, 
we show that there is a similar bound for the minimal number of {\em positive} 
self--intersection points of such an immersion, a knot invariant introduced by R. Gompf which is  called 
the positive kinkiness.

First let us recall the definition of the unknotting number and the positive respectively negative 
unknotting numbers of a link.

\begin{defn}
Let $K \subset S^3$ be an oriented link.
\begin{enumerate}
\item The unknotting number $u(L)$ of $L$ is the smallest number of crossing changes needed
to alter $L$ to the trivial link.
\item The positive (negative) unknotting number $u_+(L)$ ($u_-(L)$) is 
the smallest number of positive (negative) crossing changes in any sequence
of crossing changes deforming $L$ into the trivial link.
\end{enumerate}
\end{defn}

Here we say that a crossing change is positive if it replaces a positive crossing by a negative one, otherwise
the crossing change is called negative. Note that positive and negative unknotting numbers may really depend on the orientation 
if the link has more than one component. If a link is described as the closure of a braid, we will always assume that the orientation
is chosen such that all strings are oriented coherently, so that a generator produces a positive crossing.

In \cite{G}, R. Gompf introduced the notion of the kinkiness of a knot.
We restate his definition in a slightly modified form to include the case of a link with more than one component.
In the sequel, we will assume that 
all immersions of surfaces in the 4--ball are smooth and proper in the sense that the only singularities are transverse 
double points and that they are embeddings near the boundary.

\begin{defn} Let $L \subset S^3$ be an oriented link.
The positive (negative) kinkiness $\kappa_+(L)$ ($\kappa_-(L)$) of $L$ is the smallest number of positive 
(negative) double points of a proper immersion $F \hookrightarrow D^4$ with $\partial F=L$, where $F$ is a
connected surface of genus~0.
The kinkiness of an oriented link $L$ is the pair $(\kappa_+(L),\kappa_-(L))$.
\end{defn}

Clearly $u(L) \geq u_+(L) + u_-(L)$. Since every sequence
of unknotting operations of $p$ positive and $n$ negative crossing changes starting with $L$ and ending with the 
trivial link can be used to construct an immersion of a connected surface of genus zero in the 4--ball spanning $L$ 
which has $p$ positive and $n$ negative double points~\cite{CL}, we also have the inequalities 
$u(L) \geq \kappa_+(L)+\kappa_-(L)$ and $u_{\pm}(L) \geq \kappa_{\pm}(L)$. Observe that the kinkiness of a knot is 
clearly a concordance invariant, whereas the unknotting numbers are not.

\begin{figure}[ht]
\begin{center}
\epsfig{file=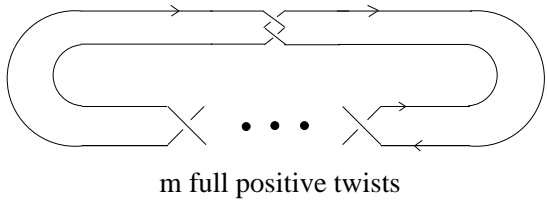}
\caption{Twist knot $K_m$}\label{twist}
\end{center}
\end{figure}

\begin{example}\label{twistknotex}
For a positive integers $m$, the twist knot $K_m$ is the knot drawn in figure~\ref{twist}. It is easy to see that
this knot can be unknotted using either $m$ negative crossing changes or one positive crossing change, therefore
$u_+(K_m)=\kappa_+(K_m)=0$ and $u_-(K_m)=\kappa_-(K_m)=0$, whereas the unknotting number is one. 
\end{example}

\begin{defn}
We will denote by
$e: B_n \rightarrow \Z$
the abelianization map sending each generator of the n--string braid group $B_n$ to $1$.
\end{defn}

In \cite{R2}, Rudolph proved the so called ``slice Bennequin inequality'', which is a lower bound for the
slice genus of a knot $K$ given as the closure of a braid $\beta \in B_n$, namely
\begin{equation}
g^*(\hat{\beta}) \geq 1 + \frac{1}{2}(e(\beta)-n-1).
\end{equation}
If the braid $\beta$ is strongly quasipositive, then this estimate is sharp.
The following result provides a similar estimate for the positive kinkiness of the closure of a braid.

\begin{thm}\label{immersedBI}
Suppose $\beta \in B_n$ is a braid of n strings. 
Let $c$ denote the number of components of its closure $\hat{\beta}$. 
Then
\begin{equation}
\kappa_+(\hat{\beta}) \geq 1 + \frac{1}{2}(e(\beta)-n-c).
\end{equation}
If moreover $F \hookrightarrow D^4$ is a proper immersion of a surface with boundary $\hat{\beta}$ which has
$p$ positive self--intersection points (and any number of negative self--intersection points) and no
closed component, then
\begin{equation}
\chi(F)-2p \leq n-e(\beta).
\end{equation}
\end{thm}

\begin{proof}
As the second equation clearly implies the first one, we only have to prove the second statement.
First let us suppose that the braid $\beta$ is positive, i.e. it can be written as a product
$\beta=\sigma_{i_1} \cdots \sigma_{i_k}$
of generators $\sigma_i$ for certain indices $1 \leq i_j \leq n-1$. Then $e(\beta)=k$. 
Let $L=\hat{\beta}$ denote the closure of $\beta$.
Since a positive braid is in particular  strongly quasipositive, 
the slice Euler characteristic $\chi_s(L)$ of $L$ equals $n-k$
by \cite{R2}. 
The link $L$ can be obtained as the transverse intersection of a complex plane curve
$f^{-1}(0) \subset \C^2$ and the 3--sphere $S^3 \subset \C^2$ \cite{R2}.
As a slight perturbation of the coefficients of $f$ does not change the link type of $L$, we can assume that
the projective closure $V \subset \CP^2$ of $f^{-1}(0)$ is a smooth algebraic curve. 
Note that $\chi_s(L)=\chi(V \cap D^4)$ by \cite{R2}. Let $V'=V \setminus int(D^4)$. Observe that $V'$ is still connected.
Now suppose that $F \hookrightarrow  D^4$ is a proper immersion of a surface having no closed component with boundary $\partial F=L$ 
which has $p$ positive double points. 
Blowing up the negative double points of $F$ and replacing the positive double points
by handles, we obtain a nullhomologous embedded surface $F' \subset D^4 \#_n \cpq$, where $n$ is the number of negative double
points of the immersion, and $\chi(F')=\chi(F)-2p$, see for instance~\cite{FS} for a detailed description of the construction. 
Now consider the union $W=V' \cup F' \subset \CP^2 \#_n \cpq$.  The Euler characteristic of $W$ is given by
\[
\begin{array}{rl}
\chi(W) = \chi(F') + \chi(V') &= \chi(F)-2p + (\chi(V) - \chi(V \cap D^4))  \\
&= \chi(F)-2p-\chi_s(L) + \chi(V).
\end{array}
\]
Now $W$ is a connected smoothly embedded surface representing the same homology class as the algebraic curve
$V$ (with respect to a suitable structure as a rational surface on $\CP^2 \#_n \cpq$),
hence $\chi(W) \leq \chi(V)$ by~\cite{KM2} or~\cite{OZ}.  
Therefore $\chi(F)-2p-\chi_s(L) \leq 0$ and the claim follows.

As to the general case, we can write $\beta$ as a product of generators and their inverse elements, i.e.
$\beta=\sigma_{i_1}^{\epsilon_1} \cdots \sigma_{i_k}^{\epsilon_k}$
with indices $1 \leq i_j \leq n-1$ and signs $\epsilon_j \in \{ -1, +1\}$. Let $a$ denote the number of
positive exponents and let $b=k-a$ denote the number of negative exponents. Then $e(\beta)=a-b$.

Consider the positive braid $\gamma=\sigma_{i_1} \cdots \sigma_{i_k}$. Clearly
$e(\gamma)=k=a+b$. The braid $\beta$ can be obtained 
from $\gamma$ by $b$ positive crossing changes, in particular $\hat{\beta}$ and $\hat{\gamma}$ have the same 
number of components. As observed in \cite{CL}, there exists a union of immersed annuli in $S^3 \times [0,1]$
connecting $\hat{\gamma}$ and $\hat{\beta}$ which has exactly $b$ self--intersection points, all of them being
positive. Now suppose that we are given an immersion $F \hookrightarrow D^4$ as above.
Gluing this immersed surface with the annuli connecting $\hat{\beta}$ and $\hat{\gamma}$
yields an immersed surface $F' \hookrightarrow D^4$ with boundary $\hat{\gamma}$ 
which has $p+b$ positive double points, note that of course $\chi(F')=\chi(F)$.
By the first part of the proof we can conclude that
\[
\chi(F)-2p-2b \leq n-a-b.
\]
As $e(\beta)=a-b$, this proves our claim.
\end{proof}

Using the inequality $u_+(L) \geq \kappa_+(L)$ for any link $L$, we immediately obtain the estimate
\begin{equation}
u_+(\hat{\beta}) \geq 1 + \frac{1}{2}(e(\beta)-n-c).
\end{equation}
a proof of which using only three--dimensional
techniques has been announced by W.W. Menasco in \cite{M}.

Now suppose we are given an oriented diagram $D$. 
A Seifert circle of $D$ is called {\em strongly negative} if 
it is not adjacent to any positive crossing. We denote by $s_-(D)$ the number of strongly negative Seifert circles
of $D$, by $s(D)$ the total number of Seifert circles and by $w(D)$ the writhe of $D$.

\begin{cor}
Let $L$ be a non--trivial oriented link with $c$ components. Then
\begin{equation}\label{estimate}
\kappa_+(L) \geq 1 + \frac{1}{2} (w(D)-s(D)-c) + s_-(D)
\end{equation}
for every diagram $D$ of $L$.
\end{cor}

\begin{proof}
If all the crossings of the diagram are negative, inspection of the Seifert surface obtained by Seifert's algorithm from $D$  shows 
that the right hand side of equation~(\ref{estimate})
is a non--positive number, so we can assume that there is at least one positive crossing. The proof is a slight modification
of an argument used in~\cite{R3} to prove a similar inequality for the slice genus. 

Let $x_\pm(D)$ denote the number of positive respectively negative crossings of $D$.
Consider the surface $Q$ which is obtained by gluing disks corresponding to those Seifert circles which are not strongly negative 
with bands corresponding
to all the positive crossings of $D$. Then $Q$ is a quasipositive surface~\cite{R3}, and its boundary $\partial Q$ is a 
strongly quasipositive link, having a positive diagram $D^+$ with $s(D)-s_-(D)$ Seifert circles and $x_+(D)$ crossings. 
By~\cite{R3}, this implies that
\[
\chi(Q)=\chi_s(\partial Q)=(s(D)-s_-(D))-x_+(D)
\]
Now suppose that $F \hookrightarrow D^4$ is a proper immersion of a connected
surface of genus zero with boundary $L$ and $\kappa_+(L)$ positive double points. Then $\chi(F)=2-c$. 
Let $S$ denote the Seifert surface for $L$ constructed by
applying Seifert's algorithm to $D$. Then $S$ is obtained from $Q$ by adding additional disks and bands, in particular $Q \subset S$.
By gluing $F$ and $S$ along $L$ and removing $Q$, we therefore obtain
a proper immersion $F' \rightarrow D^4$ with boundary $\partial Q$ which has $\kappa_+(L)$ positive double points. Moreover 
$\chi(F')=\chi(F)+\chi(S)-\chi(Q)$. Now it is known~\cite{Y} that one can deform $D^+$ into a diagram which is the closure
of a braid $\beta$ while preserving the writhe and the number of Seifert circles, i.e. $e(\beta)=w(D^+)=x_+(D)$ and $\beta$ has
$s(D)-s_-(D)$ strings. As $\hat{\beta}=\partial Q$ we obtain from Theorem~\ref{immersedBI} that
\[
\chi(F)+\chi(S)-\chi(Q)-2\kappa_+(L) \leq \chi(Q)
\]
which implies the desired inequality, note that $\chi(S)=s(D)-x_+(D)-x_-(D)$ and $w(D)=x_+(D)-x_-(D)$.
\end{proof}

\begin{cor}\label{stronglyquasipositive}
Suppose that a knot $K$ is concordant to a non--trivial strongly quasipositive knot. Then
\begin{equation}
u_+(K) \geq \kappa_+(K) \geq g^*(K) > 0.
\end{equation}
In particular the knot cannot be unknotted using only negative crossing changes.
\end{cor}

\begin{proof}
As the kinkiness and the slice genus are concordance invariants, we can asssume that $K$
itself is strongly quasipositive.
By definition, a strongly quasipositive knot $K$ can be obtained as the closure $\hat{\beta}$ of some
strongly quasipositive braid $\beta=\sigma_{i_1,j_1} \cdots \sigma_{i_k,j_k} \in B_n$
(here we use the notation from \cite{R2}). Then $e(\beta)=k$ and
as shown in \cite{R2}, the slice genus $g^*(K)$ is given
by
\[
g^*(K)=1+\frac{1}{2}(e(\beta)-n-1).
\]
Moreover the quasipositive braided Seifert surface $S(\sigma_{i_1,j_1}, \cdots ,\sigma_{i_k,j_k})$
for $\hat{\beta}$ has Euler characteristic $n-k$, so $g^*(K)=g(K) > 0$. 
Now the assertion follows from Theorem~\ref{immersedBI}.
\end{proof}

\begin{remark} In fact the arguments used in the proof of Theorem~\ref{immersedBI} show 
that the inequality $\kappa_+(K) \geq g^*(K)$
is also true for $\C$--transverse knots, i.e. for knots which can be obtained as the transverse intersection of the
3--sphere and a complex plane curve. However there are non--trivial $\C$--transverse knots which are slice,
see Example 3.2 in \cite{R1}.
\end{remark}

\begin{remark}  Note that the slice genus and the positive kinkiness are concordance invariants whereas the property
of being strongly quasipositive is not. In fact, suppose that $K$ is some non--trivial strongly quasipositive knot. Pick a slice
knot $K_0$ with genus $1$ and consider the knot $K'=K \# K_0$. Then $K$ and $K'$ represent the same concordance class, but the genus
of $K'$ is $g(K)+1$. Therefore the slice genus of $K'$ -- which is $g(K)$ -- is strictly less than its genus, which, by the
results in~\cite{R2}, implies that $K'$ is not strongly quasipositive.
\end{remark}

\begin{cor}\label{positive} 
Suppose that a knot $K$ is the closure of a positive braid. Then
the kinkiness of $K$ is $(g(K),0)$ and $u(K)=u_+(K)=g(K)$.
\end{cor}

\begin{proof}
Suppose $K$ is the closure of the positive braid $\beta \in B_n$. Let $k$ denote the number of letters in $\beta$.
It is well known that $k$ can be unknotted using $\frac{1}{2}(k-n+1)$ crossing changes, see for instance
\cite{BW}. However, since $\beta$ is positive, $e(\beta)=k$, and therefore we can conclude, using once more the results
from \cite{R2}, that
\[
u(K) \leq \frac{1}{2}(e(\beta)-n+1)=g^*(K)=g(K).
\]
Since it is always true that $u(K) \geq g^*(K)$, these two numbers must coincide.
Now, by Corollary~\ref{stronglyquasipositive}, we also have $\kappa_+(K) \geq g(K)$. Hence
\[
u(K) \geq u_+(K) \geq \kappa_+(K) \geq g(K)=u(K)
\]
and we obtain that also $u_+(K)=\kappa_+(K)=g(K)$. Finally the inequality $u(K) \geq u_+(K)+u_-(K)$
implies that $u_-(K)=0$, and as $u_-(K) \geq \kappa_-(K)$ we have $\kappa_-(K)=0$, as claimed.
\end{proof}

\begin{example}
Suppose that $p$ and $q$ are coprime positive integers. The $(p,q)$--torus knot $t(p,q)$ is the closure of the
$q$--string braid $(\sigma_1\sigma_2 \cdots \sigma_{q-1})^p$, in particular it is the closure of a positive braid.
It is well known that the unknotting number of $t(p,q)$ equals its genus, this follows from the positive solution of 
the Milnor--conjecture by Kronheimer and Mrowka~\cite{KM1}. By Corollary~\ref{positive}, the genus also equals the 
positive unknotting number $u_+$ and the positive kinkiness, in particular every sequence of crossing changes turning 
$t(p,q)$ into the unknot which has minimal length consists of positive crossing changes only.
\end{example}

\begin{figure}[ht]
\begin{center}
\epsfig{file=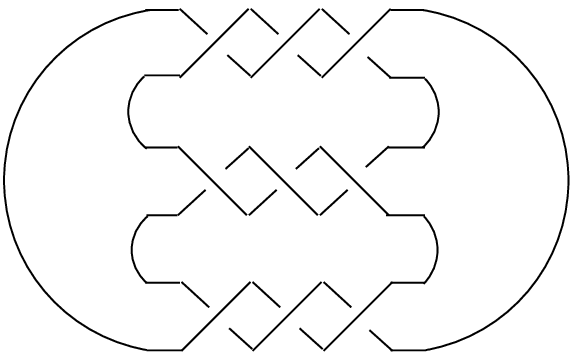}
\caption{K(-3,3,-3)}\label{pretzel}
\end{center}
\end{figure}

At this point is seems necessary to fix a convention as to the notation for pretzel knots.
Suppose that $p,q$ and $r$ are odd integers. As in~\cite{K} and~\cite{CL},
we denote by $K(p,q,r)$ the pretzel
knot which is the boundary of the surface obtained by gluing two disks with 3 vertical bands, where the first band
has $|p|$ twists, positive if $p>0$ and negative if $p<0$, similarly the second and third band have $|q|$
resp. $|r|$ twist. The pretzel knot $K(-3,3,-3)$ is drawn in Figure~\ref{pretzel}.

Note that, with this convention, a positive $p$ gives rise to $p$ negative crossings, so for instance the
knot $K(1,1,1)$ has three negative crossings, in fact it is the left--handed trefoil knot.

\begin{example}
Corollary~\ref{positive} is in general not true for strongly quasipositive knots. For an example, consider the pretzel knot 
$K=K(-9,5,-9)$. It has been proved in~\cite{R2} that 
this knot is strongly quasipositive (the notation for pretzel knots in~\cite{R2} differs from
our notation, so $K$ is denoted by $P(9,-5,9)$ in~\cite{R2}), in particular
it is not slice. Also note that the signature of $K$ is zero, so it is not
positive by~\cite{P}, in particular it cannot be the closure of a positive braid.
As the genus of $K$ is clearly one we obtain that $g^*(K)=1$. However, as shown in Example 2.12 in \cite{CL},
the unknotting number of $K$ is at least two. 
\end{example}

In~\cite{CL}, Cochran and Lickorish  showed that, as a consequence of Donaldson's Theorem on the intersection forms of smooth 4--manifolds, 
certain pretzel knots cannot be unknotted using only positive crossing changes. Using
Corollary~\ref{stronglyquasipositive}, we are now able to prove that this holds for a larger class of pretzel knots.

\begin{cor}\label{pretzelknots}
Assume that $p,q,r$ are odd integers and $\min \{ p+q,q+r,p+r\} > 0$. Then the pretzel knot $K(p,q,r)$ has
infinite order in the smooth knot concordance group and
cannot be unknotted using only positive crossing changes. 
\end{cor}

\begin{proof} The first part of the statement is an immediate consequence of the results
in~\cite{R2}. As shown there, the mirror image of such a pretzel knot is strongly quasipositive (again note that
in \cite{R2}, a slightly different definition of pretzel knots is used, in fact the knot which is there called
the pretzel knot of type $(p,q,r)$ is the mirror image of the knot which we denote by $K(p,q,r)$).
The sum of
two strongly quasipositive knots is again strongly quasipositive,
and hence every non--trivial strongly quasipositive knot has infinite order in the smooth concordance group by~\cite{R2}.
Using Corollary~\ref{stronglyquasipositive}, we can also conclude that the negative kinkiness $\kappa_-(K(p,q,r))$ is
not zero and that $K(p,q,r)$ cannot be unknotted by using only positive crossing changes.
\end{proof}

\begin{remark}
Note that, by combining the arguments given in~\cite{CL} with~\cite{FS},  Corollary~\ref{pretzelknots} 
can also be proved without referring to the theory of quasipositive knots. 

In fact, as explained in \cite{CL}, a pretzel knot of the form $K(-1-2n,-1+2n,-1-2n)$ for some non--negative
integer $n$ is the boundary of an embedded disk $D$ in $\CP^2\#_n\cpq \setminus D^4$ whose homology class
is $(3,1, \cdots, 1)$. Now suppose $K(-1-2n,-1+2n,-1-2n)$ were the boundary of an immersed disk in the
4--ball having only negative double points. By gluing this disk with $D$ we could then construct an
immersed sphere in the rational surface $\CP^2 \#_n\cpq$ representing the class $(3,1, \cdots, 1)$ which has
only negative double points, in contradiction to the results in \cite{FS}. Hence we can conclude that a pretzel knot
of the form $K(-1-2n,-1+2n,-1-2n)$ has non--zero positive kinkiness. Passing to the mirror image we obtain that
the knots $K(1+2n,1-2n,1+2n)$ have non--zero negative kinkiness. The same then holds for all pretzel knots of
the form $K(1+2n+a,1-2n+b,1+2n+c)$ for non--negative integers $a,b,c,n$, because such a knot can be obtained from
$K(1+2n,1-2n,1+2n)$ by positive crossing changes. In particular, none of these knots can be unknotted
using only negative crossing changes, this result has also been obtained
in~\cite{CL} in the case that $n \leq 8$.

Now suppose we are given three odd integers $p,q,r$
such that $\min \{p+q,p+r,q+r\} > 0$. Then at least two of these numbers, say $p$ and $r$, are positive, and we can
also assume that $p \leq r$. If we define integers $b,c,n$ by $p=1+2n$, $r=p+c$ and $q=2-p+b$, then $b,c$ and
$n$ are non--negative and $K(p,q,r)=K(1+2n,1-2n+b,1+2n+c)$. So $K(p,q,r)$ has non--zero negative kinkiness, hence 
it cannot be unknotted by negative crossing changes only. As it can clearly be unknotted by using only positive
crossing changes we can conclude that its positive kinkiness is zero.
By Theorem 1.1 in \cite{G}, this implies
that $K(p,q,r)$ has infinite order in the knot concordance group.
Note that this proof, as the proof via Theorem~\ref{immersedBI}, is based on the positive solution of the Thom conjecture 
respectively its immersed version.
\end{remark}

\begin{example}
Suppose that $p,q,$ and $r$ are three odd integers such that $pq+pr+qr=-1$ and consider the pretzel knot $K(p,q,r)$.
The Alexander polynomial of $K(p,q,r)$ is then trivial and the knot is topologically slice by \cite{F}.
First assume that two of the numbers $p,q,r$ are positive. It is then not hard to see that either 
$\{-1,1\} \subset \{p,q,r\}$ which implies that $K(p,q,r)$ is trivial, or $\min \{p+q,p+r,q+r\} > 0$. By 
Corollary~\ref{pretzelknots} we can conclude that either $K(p,q,r)$ has infinite order in the knot concordance group
and cannot be unknotted using only positive crossing changes or $K(p,q,r)$ is trivial, 
thus recovering Corollary 4.2 and Corollary 4.3 in \cite{CL}.
If two of the numbers $p,q$ and $r$ are negative the same conclusion holds for the mirror image $K(-p,-q,-r)$.
\end{example}

In some cases Corollary~\ref{stronglyquasipositive} can be used to
show that certain knot concordance classes do not contain any strongly quasipositive knots. In fact one can find
a subgroup of infinite rank in the smooth knot concordance group in which all concordance classes have this property.

\begin{cor}
There exists a subgroup of infinite rank in the smooth knot concordance group of which no non--trivial element can be
represented by a strongly quasipositive knot.
\end{cor}

\begin{proof}
Consider the subgroup $G$ of the smooth knot concordance group generated by all twist knots $K_m$, $m>0$.
Recall (see Example~\ref{twistknotex}) that $\kappa_+(K_m)=\kappa_-(K_m)=0$. This is clearly also true for every knot
which can be written as a connected sum of twist knots and mirror images thereof. As the kinkiness is a concordance invariant,
this implies that for every knot $K$ whose concordance class is in $G$, $\kappa_-(K)=\kappa_+(K)=0$. If such a knot is
strongly quasipositive, then Corollary~\ref{stronglyquasipositive} shows that $g^*(K)=0$, in other words $K$ is slice
and the concordance class represented by it is the trivial one. In~\cite{J}, B.J. Jiang modified the
arguments of~\cite{CG} to show that 
$G$ has infinite rank, this is even true for the subgroup of $G$ generated by all twist knots $K_m$ where $4m+1$ is a
square. Hence $G$ has all the required properties and the proof is complete.
\end{proof}

\end{document}